# STATIONARY ITERATION METHODS FOR SOLVING 3D ELECTROMAGNETIC SCATTERING PROBLEMS

ALEXANDER SAMOKHIN[*], YURY SHESTOPALOV[†], AND KAZUYA KOBAYASHI[‡]

**Abstract.** Generalized Chebyshev iteration (GCI) applied for solving linear equations with nonselfadjoint operators is considered. Sufficient conditions providing the convergence of iterations imposed on the domain of localization of the spectrum on the complex plane are obtained. A minimax problem for the determination of optimal complex iteration parameters is formulated. An algorithm of finding an optimal iteration parameter in the case of arbitrary location of the operator spectrum on the complex plane is constructed for the generalized simple iteration method. The results are applied to numerical solution of volume singular integral equations (VSIEs) associated with the problems of the mathematical theory of wave diffraction by 3D dielectric bodies. In particular, the domain of the spectrum location is described explicitly for low-frequency scattering problems and in the general case. The obtained results are discussed and recommendations concerning their applications are given.

**1. Introduction.** Analysis of the scattering of electromagnetic waves by three-dimensional inhomogeneous anisotropic dielectric structures is of crucial importance for studying various applied and theoretical issues. In the short-wave range, the asymptotical methods have proved to be the most efficient technique. In the resonance and low-frequency ranges, it is necessary to solve exact equations that describe the wave scattering; namely, the Maxwell equations with radiation conditions at infinity or volume singular integral equations (VSIEs) with respect to the wave field in the dielectric structure. Partial differential equations (PDEs) may seem to be more appropriate for numerical solution because their discretization results in a system of linear algebraic equations (SLAE) with a sparse matrix, while in the case of a VSIE this matrix is dense. However, when a scattering problem is considered, the solution must satisfy the condition at infinity; therefore, in order to provide the required accuracy the unknown wave field must be determined in a domain which is much larger than the actual dielectric scatterer. Finally, taking into account that the scattering problem is three-dimensional, we obtain a SLAE of extremely high size when PDEs are discretized and solved numerically by a finite-difference (FD) or a finite-element method (FEM). Application of approximate conditions at infinity often leads to significant loss of accuracy that cannot be controlled.

In view of these difficulties, we develop and apply the VSIE method. Using the fast discrete Fourier transform and taking into account that the VSIE kernels depend on the difference of arguments, we consider fast algorithms of the matrix-vector multiplication. Next, using iteration techniques, we perform efficient numerical solution of the initial scattering problems by the VSIE method [9].

Recall main parameters that govern efficiency of a numerical algorithm: number $T$ of arithmetic operations which is necessary for obtaining the sought-for solution with the required accuracy and volume $M$ of the memory required for implementation of the algorithm.

Matrix-vector multiplication is the most laborious operation when iteration techniques are applied. Therefore, the number of these multiplications in the process of implementation of the algorithm will be called *the number of iterations*. The value of $T$ may be estimated from the formula
$$T \approx L(T_A + T_0) + T_M.$$
Here $L$ is the number of iterations, $T_A$ is the number of arithmetic operations for multiplying a matrix by a vector, $T_0$ is the number of arithmetic operations for other calculations performed within one iteration, and $T_M$ is the number of arithmetic operations for forming the SLAE matrix when an integral equation is discretized. As a rule, $T_0 << T_A$ и $T_M << L T_A$.

---
[*]Moscow State Technical University of Radio Engineering and Automation, Vernadsky av. 78, 117648 Moscow, Russia (absamokhin@yandex.ru)
[†]Karlstad University, SE-651 88Karlstad, Sweden (youri.shestopalov@kau.se)
[‡]Chuo University 1-13-27 Kasuga, Bunkyo-ku, Tokyo 112-8551 Japan (kazuya@tamacc.chuo-u.ac.jp)

The memory volume required for implementation of the algorithm can be estimated by the expression

$$M \approx M_A + M_{ITER}$$

Here $M_A$ and $M_{ITER}$ are the memory volumes required, respectively, for storing the SLAE matrix or the corresponding array and implementing the iteration procedure.

Quantities $T_A$, $T_M$, and $M_A$ depend solely on the method of dicretization of the integral equation, while $L$, $T_0$, and $M_{ITER}$ are governed by the iteration algorithm in use.

Since the VSIE kernels depend on the difference of arguments, we can provide the fulfillment of the condition $M_A \sim N$, where *N* is the matrix size using efficient dicretization.

The value of $T_A$ is another important characteristic of the algorithm. If no special algorithms are used, then $T_A \approx N^2$, which causes abnormal computational expenses due to huge matrix dimensions. Using the fast discrete Fourier transform, it is possible to create an efficient algorithm with $T_A \sim N\,LOG(N)$, where $LOG(N)$ is a function taking values on the set of whole numbers and equal to the sum of all prime factors of number *N*. Obviously, this is a slowly varying function with respect to increasing; for example, $LOG(10^6) = 42$ and $LOG(10^8) = 56$.

Now consider minimization of the number of iterations *L* and quantities $T_0$ and $M_{ITER}$ which are governed by the iteration algorithm. Usually, the integral equations under study are solved by the iteration algorithm called *Generalized Minimal Residual Algorithm* (GMRES) and its various modifications [4]. When this method is applied, the iterations parameters are determined in the course of calculations and depend on the current iteration index. Such methods may be called nonstationary iteration algorithms. GMRES-type algorithms are very popular; however, their implementation requires comparatively large amount of computer memory with $M_{ITER} \approx n\,N$, where *n* is the dimension of the Krylov subspace (note that the rate of convergence of iterations to the solution increases together with *n*). Also, the number of iterations required for providing reasonable accuracy of solution is often very large.

There is another family of iteration methods called stationary iteration algorithms for which the iterations parameters are determined before the iteration procedure. For a stationary iteration algorithm, the values of $T_0$ и $M_{ITER}$ turn out to be minimal. These methods are usually applied for solving equations with selfadjoint operators for which the boundaries of the spectrum location on the real axis of the complex plane are known. The Chebyshev iteration and the method of simple iteration [6] belong to this family. For this methods we have $M_{ITER} \approx N$.

Iteration techniques, mainly GMRES and its modifications, are widely used for numerical solution of VSIEs [ 6, 7,11,13]. Our goal is to work out efficient iteration methods for numerical solution of VSIEs in the *low-frequency range*. The problems of this class find broad applications in various domains; however, unlike high-frequency range, the specific features connected with this particular range of frequencies have not been taken into account as far as efficient numerical approaches are concerned. The Born iteration method constitutes here the only exception. However, this method can be used solely when the permittivity of the dielectric scatterer is close to that of free space; this condition is a very severe limitation which does not hold in the majority of applied problems.  In this work, we propose to fill this gap and apply the generalized Chebyshev iterations and the method of simple iteration to solve VSIEs numerically in the low-frequency range. The latter requires the knowledge concerning the location of the VSIE operator spectrum on the complex plane, which in its turn may be considered as a urgent problem of algebra and numerical analysis. In fact, all modern iteration techniques aimed at solution of very big linear equation systems employ the information on the operator spectrum which enables one to improve significantly the algorithm productivity. We describe *explicitly* the domain of the spectrum location on the complex plane for low-frequency VSIEs *without any restrictions* imposed on the permittivity function. This result plays a crucial role for creating  efficient methods and algorithms for numerical solution of VSIE and the corresponding scattering and diffraction problems. The very recent findings [7, 14, 15] also show that they are extremely important for the solution to inverse problems of reconstructing permittivity of 3D dielectric bodies.

It should be noted also that utilization of modern multi-processor computers leads to substantial increase in the rate and performance of computations when VSIE is applied because that latter admits excellent parallelization for very large and dense matrices.

In Section 2 we describe the generalized method of simple iteration (GSI) for solving equations with nonselfadjoint operators [9] and construct a finite algorithm of finding an optimal iteration parameter in the case of arbitrary location of the operator spectrum on the complex plane. In Section 3 we describe the generalized Chebyshev iteration (GCI) applied for solving equations with nonselfadjoint operators. These techniques demand the knowledge of the spectrum localization on the complex plane. It is not possible to obtain this information for the majority of problems. However, one can do it in a number of important particular cases. In Section 4 we formulate the problems of electromagnetic scattering and reduce them to VSIEs. In Section 5 we generalize the results obtained in [1] and describe the domain of the spectrum location explicitly for low-frequency scattering problems and in the general case. In Section 6 we illustrate theoretical findings of this study by the results of computations obtained for a representative scattering problem, namely, low-frequency scattering of a plane electromagnetic wave from an inhomogeneous dielectric ball, and compare our iterative strategies and methods employing optimal parameters with some typical results obtained by a state-of-the-art algorithm. In Section 7 we discuss the obtained results and give recommendations concerning their applications.

**2. Generalized simple iterations.** Consider a linear operator equation

(2.1) $$\hat{A}u = f,$$

in the Banach space $E$, where $\hat{A}$ is a bounded and therefore continuous operator (generally, nonselfadjoint).

Write equation (2.1) in the equivalent form

(2.2) $$u - \hat{B}_\mu u = f/\mu.$$

Here $\hat{B}_\mu$ is a linear operator defined by

(2.3) $$\hat{B}_\mu = \frac{\mu \hat{I} - \hat{A}}{\mu},$$

and $\mu \neq 0$ is an arbitrary complex number.

DEFINITION. *The number*

$$\rho_0 = \sup|\eta|, \quad \eta \in \sigma(\hat{B}),$$

*is called the spectral radius of operator* $\hat{B}$. *Here* $\sigma(\hat{B})$ *the spectral set of* $\hat{B}$ *on the complex plane; that is, a set of points* $\eta$ *such that operator* $(\hat{B} - \eta \hat{I})$ *has no inverse defined on the whole space*

The following statement is proved in the theory of linear operators [5].

THEOREM 2.1. *The linear operator equation* $u - \hat{B}u = f$ *has the unique solution for every* $f \in E$ *and the successive approximations*

(2.4) $$u_{n+1} = \hat{B}u_n + f, \quad n = 0,1,...$$

*converge to the solution at any initial approximation* $u_0 \in E$ *if the spectral radius of the operator satisfies* $\rho_0(\hat{B}) < 1$. *Conversely, if iterations (2.4) converge at any* $u_0, f \in E$, *then* $\rho_0(\hat{B}) \leq 1$. *The following estimate is valid for the rate of convergence*

(2.5) $$\|u_n - u\| \leq C(\rho_0(\hat{B}))^n, \quad C = const.$$

Note that if $\hat{B}$ is a normal operator then one can set $C = 1$ in (2.5).

Theorem 2.1 yields the convergence of the successive approximations

(2.6) $$u_{n+1} = \hat{B}_\mu u_n + f/\mu, \quad n = 0,1,...$$

to the solution of (2.2) and hence of (2.1) at any $u_0, f \in E$ if the spectral radius of operator $\hat{B}_\mu$ is less than one, that is,

(2.7) $$\rho_0(\mu) = \sup |\eta(\mu)| < 1, \quad \eta(\mu) \in \sigma(\hat{B}_\mu).$$

From (2.3) it follows that iterations (2.6) can be represented in the form

(2.8) $$u_{n+1} = u_n - \frac{1}{\mu}(\hat{A}u_n - f), \quad n = 0,1,...$$

Our next goal will be to answer the question: for which location of the spectrum of initial operator $\hat{A}$ on the complex plane one can determine complex numbers $\mu$ such that iterations (2.8) converge to the solution of (2.1).

It is easy to show that there is a one-to-one correspondence between the spectra of operators $\hat{A}$ and $\hat{B}_\mu$ given by the expression

(2.9) $$\eta = (\mu - \lambda)/\mu, \quad \lambda \in \sigma(\hat{A}), \quad \eta \in \sigma(\hat{B}_\mu).$$

Using (2.7) and (2.9) and geometric considerations one can prove the following [9]

THEOREM 2.2. *The complex numbers $\mu$ at which iterations (2.8) converge to the solution of (2.1) at any $u_0, f \in E$, exist if and only if the origin of the complex plane is situated outside the convex envelope of the spectrum of operator $\hat{A}$.*

Here we used the following definition.

DEFINITION. *The convex envelope of M is intersection of all convex sets containing M.*

Assume that the condition of Theorem 2.2 holds. Then a natural question arises concerning the determination of iteration parameter $\mu_0$ which provides the best convergence rate. From (2.7), (2.9), and Theorem 2.1 it follows that iterations (2.8) would converge to the solution with the rate of geometrical progression having the common ratio

(2.10) $$\rho_0(\mu) = \frac{\sup|\mu - \lambda|}{|\mu|}, \quad \lambda \in \sigma(\hat{A}).$$

The best convergence rate will be when $\mu$ takes the value at which function $\rho_0(\mu)$ attains minimum. Let $S_\mu$ denote a circle on the complex plane with the origin at $\mu$ with the minimal radius $R = \sup|\mu - \lambda|, \lambda \in \sigma(\hat{A})$ that envelops all points of the spectrum of operator $\hat{A}$. Draw the tangent lines to $S_\mu$ from the origin and denote by $\alpha$ the angle between them (see Fig. 2.1 where the spectrum occupies the shaded region). Then formula (2.10) yields $\rho_0(\mu) = \sin(\alpha/2)$. Thus, we have proved the following statement.

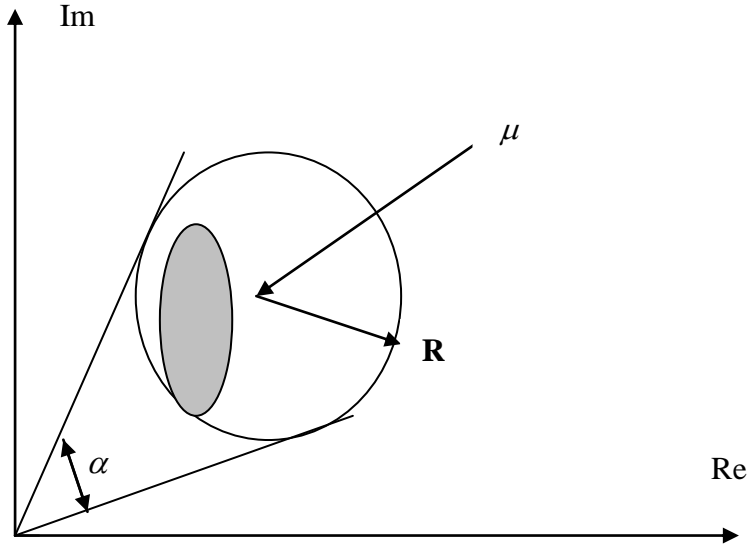

**Figure 2.1:** The choice of the optimal iteration parameter.

THEOREM 2.3. *Let the origin of the complex plane be situated outside the convex envelope of the spectrum of operator $\hat{A}$. Denote by $S_0$ a circle on the complex plane which contains all points of the spectrum of operator $\hat{A}$ and "can be seen" from the origin at a minimal angle $\alpha_0$. Then the best convergence rate of iterations (2.8) to the solution of equation (2.1) is attained when $\mu_0$, coincides with the origin of circle $S_0$. The iterations will converge to the solution with the rate of geometrical progression having the common ratio $\rho_0 = \sin(\alpha_0/2)$.*

Note that if $\hat{A}$ is a selfadjoint and positive definite operator, then these results are well-known. Indeed, in this case the spectrum of the operator is situated on the positive real semi-axis on the complex plane which yields, by virtue of Theorem 2.3, $\mu_0 = (M+m)/2$, where $M$ and $m$ are the upper and lower bounds of the spectrum, and coincides with the classical result.

All the proofs for the iteration method under consideration are geometric since the use of purely analytical techniques for nonselfadjoint operators is very complicated. We apply the same approach for constructing a finite algorithm for determination of the sought-for circle $S_0$ and the corresponding iteration parameter $\mu_0$ for a given convex envelope of the spectrum which does not contain the origin. Assume to this end without loss of generality that this envelope is an arbitrary convex polygon with $n$ vertices on the complex plane.

The following propositions are valid.

PROPOSITION 2.1. *Assume that the spectrum of the operator is localized on a rectilinear segment. Then $\mu_0$ will be a point of intersection of the medial vertical to the segment and a circle drawn through its endpoints and the origin so that the segment will be a chord of the sought-for circle $S_0$.*

PROPOSITION 2.2. *For a polygon with n vertices the sought-for circle $S_0$ will go at least through two vertices of the polygon.*

*Proof of Proposition 2.1.* Clearly, this segment will be a chord of the sought-for circle. Next, consider minimum of $\sin(\alpha/2)$ with respect to point $x$ situated on the medial vertical taking into account that the endpoints belong to the boundary of the circle. If the segment lies on a beam drawn through the origin then $\sin(\alpha/2)$ has one extremum which is a point of minimum, namely, the medial point. If the segment does not lie on these beams, then $\sin(\alpha/2)$ has two extremum points coinciding with the points of intersection of the medial vertical to the segment and the circle that goes through the endpoints and the origin. A circle that does not contain the origin corresponds to the minimum which proves the proposition (see Fig. 2.2).

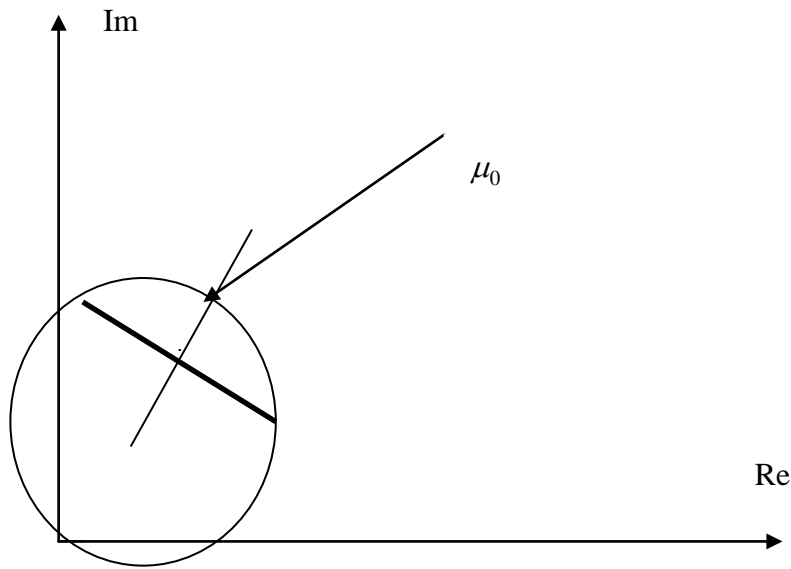

**Figure 2.2:** Optimal iteration parameter for a segment

*Proof of Proposition 2.2.* Clearly, at least one point of the polygon will belong to the boundary of the sought-for circle; otherwise, its radius may be decreased without shifting the origin. Assume that not more than one vertex lies on the circumference and join the origin of the circle with this vertex of the polygon (origin of the circle $\mu_1$ in Fig. 2.3). Then, by shifting the origin towards the vertex one can construct a smaller circle contained in the initially constructed circle and containing the polygon (origin of the circle $\mu_2$ in Fig. 2.3). The latter constitutes a contradiction which proves the proposition.

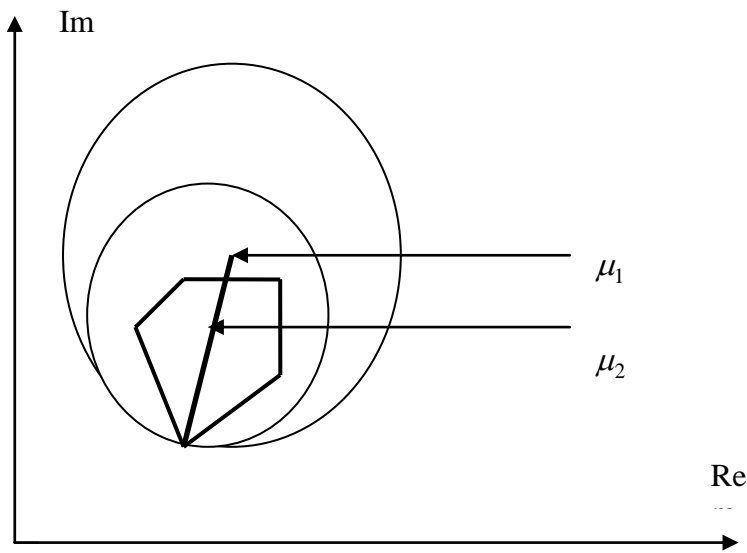

**Figure 2.3:** Illustration to the proof of Proposition 2.2.

Describe a finite algorithm of finding, for a given convex polygon, the sought-for circle $S_0$ and therefore the best iteration parameter $\mu_0$.

Step A. *The number of all segments joining the vertices of the polygon is n(n-1)/2. We consider successively these segments and construct the best circle for each of them according to Proposition 1. If for a certain segment this circle contains the whole polygon then this circle will be the sought-for $S_0$ and the algorithm stops. Note that there may be several such segments, however, the constructed circles will be the same. If none of the circles contains the whole polygon we go to step B.*

Step B. *In this case, the sought-for circle $S_0$ will go through three vertices of the polygon. The number of all triangles constructed on vertices of the polygon is n(n-1)(n-2)/6. The sought-for $S_0$ will be a circle inscribed in one or several triangles. Note that there may be several such triangles, however, the constructed circles will be the same. The algorithm terminates.*

*Proof of step A of the algorithm.* Let $\sigma_1$ and $\sigma_2$ be two sets of points on the complex plane whose convex envelopes do not contain the origin and $\sigma_1 \subset \sigma_2$. Obviously, $S_0(\sigma_1) \subseteq S_0(\sigma_2)$ and $\alpha_0(\sigma_1) \leq \alpha_0(\sigma_2)$, where $S_0(\sigma)$ is a circle containing set $\sigma$ and can be seen from the origin at a minimal angle $\alpha_0(\sigma)$. Denote by $\sigma_1$ a segment of the polygon and by $\sigma_2$ the whole polygon. It is clear that if for a certain segment the constructed circle contains the whole polygon, then this circle will be the sought-for $S_0$.

*Proof of step B of the algorithm.* Assume that the sought-for circle $S_0$ has only two vertices $\alpha$ and $\beta$ of the polygon on the boundary. Then, since the algorithm has not terminated on step A, the origin of this circle is shifted along the medial vertical with respect to the origin $\mu_{\alpha\beta}$ of best circle for a segment joining points $\alpha$ and $\beta$. By shifting the origin towards $\mu_{\alpha\beta}$, one can construct a circle that can be seen from the origin at a minimal angle and contains the whole polygon. Thus, we have arrived at a contradiction which means that the sought-for $S_0$ has at least three vertices of the polygon on the circumference.

Obviously, the constructed circle is unique.

**3. Generalized Chebyshev iteration.** GSI employs only one iteration parameter. However, for solving an equation with a selfadjoint operator the Chebyshev iteration turns out to be more efficient. In this method, the iteration parameters depend on the iteration number. Below, we apply this method for solving equations with nonselfadjoint operators. Consider a linear operator equation (2.1) in the Banach space $E$ assuming to this end that $\hat{A}$ is a bounded operator and there exists a bounded inverse $\hat{A}^{-1}$.

We seek a solution to equation (2.1) using an iteration procedure

$$u_{m+1}^l = u_m^l - \tau_{m+1}(\hat{A}u_m^l - f), \quad m = 0,1,...,n-1, \ l = 0,1,...$$

(3.1) $$u_0^0 = u_{(0)}, \quad u_0^{l+1} = u_n^l, l \geq 1,$$

where $u_{(0)}$ is an initial approximation and $\tau_m$ are complex iteration parameters. From (3.1) we see that the iteration procedure constitutes a sequence of layers and each layer consists of $n$ iterations with an identical set of parameters.

Denote by $h_m^l$ the iteration residuals defined according to the formulas

(3.2) $$h_m^l = (\hat{A}u_m^l - f), \quad m = 0,1,...n-1, \ l = 0,1,...$$

Now (3.1) and (3.2) yield

(3.3) $$h_{m+1}^l = h_m^l - \tau_{m+1}\hat{A}h_m^l, m = 0,1,...n-1, \ l = 0,1,...$$

From (3.3) it follows that within each $l$-th layer the first and the last residuals are coupled by a relationship

$$(3.4) \qquad h_n^l = \prod_{m=1}^{n}(\hat{I} - \tau_m \hat{A})\, h_0^l.$$

Successively applying (3.4) we obtain that after completing $k$ layers (3.1) residual $h_n^k$ is coupled with initial residual $h_0^0$ by the inequality

$$(3.5) \qquad \|h_n^k\| \le \left\|[\prod_{m=1}^{n}(\hat{I} - \tau_m \hat{A})]^k\right\| \|h_0^0\|.$$

Estimate (3.5) implies that the iteration parameters must minimize the norm $\left\|[\prod_{m=1}^{n}(\hat{I} - \tau_m \hat{A})]^k\right\|$.

Denote

$$(3.6) \qquad \hat{B}_k(\tau_1,...,\tau_n) = [\prod_{m=1}^{n}(\hat{I} - \tau_m \hat{A})]^k.$$

For any bounded operator $\hat{B}$ the following estimate [5] holds

$$(3.7) \qquad \|\hat{B}\| \le C\rho_0(\hat{B}),\ \ C = const,$$

where $\rho_0(\hat{B})$ is the spectral radius of operator $\hat{B}$. If $\hat{B}$ is a normal operator then $C=1$.

Applying (3.6), (3.7), and the theorem of the spectrum mapping [3], we obtain

$$(3.8) \qquad \|\hat{B}_k(\tau_1,...,\tau_n)\| \le C(\max_{z \in \sigma(A)}[\left|\prod_{m=1}^{n}(1 - \tau_m z)\right|]^k).$$

The following obvious inequalities are valid

$$(3.9) \qquad \|u - u_n^k\| = \|\hat{A}^{-1}\hat{A}(u - u_n^k)\| \le \|\hat{A}^{-1}\|\|f - \hat{A}u_n^k\| = \|\hat{A}^{-1}\|\|h_n^k\|,$$

where $u$ is the sought-for solution to equation (2.1). Now from (3.5), (3.6), (3.8), and (3.9) it follows that

$$(3.10) \qquad \|u - u_n^k\| \le \|\hat{A}^{-1}\| C\, ([\max_{z \in \sigma(A)}\left|\prod_{m=1}^{n}(1 - \tau_m z)\right|]^k)\|h_0^0\|$$

Therefore, the determination of an optimal set of iteration parameters $\tau_1,...,\tau_n$ reduces to determination of a complex polynomial of degree $n$

$$(3.11) \qquad P_n(z) = \prod_{m=1}^{n}(1 - \tau_m z),$$

which has a minimal maximum of the absolute value in the domain of location of the spectrum of operator $A$ on the complex plane. Namely, we have to solve a minimax problem

$$(3.12) \qquad \rho_{opt} = \min_{\{\tau_m\}} \max_{z \in \sigma(A)}\left|\prod_{m=1}^{n}(1 - \tau_m z)\right|$$

with respect to $\{\tau_m\}$. Obviously, the iteration parameters do not vanish. Therefore, performing the change of variables $\tau_m = 1/\mu_m$, we can represent problem (3.12) in the form

$$(3.13) \qquad \rho_{opt} = \min_{\{\mu_m\}}\left[\frac{1}{\left|\prod_{m=1}^{n}\mu_m\right|} \max_{z \in \sigma(A)}\left|\prod_{m=1}^{n}(\mu_m - z)\right|\right]$$

From (3.13) we see that the case $n = 1$ corresponds to GSI; the determination of the corresponding iteration parameter was considered above.

According to (3.10) iteration procedure (3.1), (3.12) converges if $\rho_{opt} < 1$. If the convex envelope of operator $\hat{A}$ does not contain the origin of the complex plane, then, taking into consideration Theorem 2.3, we have $\rho_{opt} \leq [\rho_0(\mu_0)]^n < 1$, where $\mu_0$ is the GSI iteration parameter. Thus we have proved the following statement.

THEOREM 3.1. *If the convex envelope of operator $\hat{A}$ does not contain the origin of the complex plane, then there exists a set of complex iteration parameters $\tau_1,...,\tau_n$ defined by (3.12) at which $\rho_{opt} < 1$ and iteration procedure (3.1) converges to the solution of equation (2.1).*

Note that unlike Theorem 2.2, Theorem 3.1 constitutes only sufficient conditions providing convergence of GCI.

In general, for the arbitrary spectrum localization on the complex plane the solution to problem (3.12) or (3.13) is not known if $n > 1$.

Consider several particular cases.

Assume that the spectrum of the operator lies on a segment of the beam drawn from the origin of the complex plane. If the spectrum is located on a segment $[a, b]$, $b>a>0$, of the real axis, then we arrive at a well-known classical solution: the sought-for real iteration parameters $\tau_1,...,\tau_n$ are the roots of the Chebyshev polynom. Then the real iteration parameters $\mu_1,...,\mu_n$ are determined from the formula [6]

$$(3.14) \qquad \mu_m = \frac{(b+a)}{2} + \frac{(b-a)}{2}\cos\frac{(2m-1)\pi}{2n}, \quad m=1,...,n.$$

We demonstrate expressions for parameters $\mu_m$ because for them the formulas are algebraically more convenient than for $\tau_m = 1/\mu_m$. From (3.14) it is clear that parameters $\mu_m$, $m=1,...,n$ belong to segment $[a, b]$.

Let the spectrum of the operator be located on a segment of a beam drawn from the origin of the complex plane. Then formula (3.13) yields $\mu_1 \exp(i\varphi),...,\mu_n \exp(i\varphi)$, where $\varphi$ is the angle of inclination of the beam to the real axis of the complex plane and $\mu_m$ are given by (3.14).

Assume now that the spectrum of the operator is contained in a circle $S_0$ and the origin of the complex plane is situated outside this domain. We prove that in this case, all iteration parameters are the same for any *n*:

$$(3.15) \qquad \tau_1 = ... = \tau_n = 1/\mu_0,$$

where $\mu_0$ is the origin of circle $S_0$ on the complex plane.

Draw a beam from the origin of the complex plane through the origin of the circle and denote by $D_1$ the diameter on the beam. By $D_2$ denote the diameter of the circle perpendicular to $D_1$. Consider the case when the spectrum of the operator is located in domain $D_1$. Then the iteration parameters $\mu_m$, $m=1$,..., *n*, lie, by virtue of (3.14), on segment $D_1$. Consider two very narrow ellipses constructed on diameters $D_1$ and $D_2$ denoting them by $\tilde{D}_1$ and $\tilde{D}_2$. Consider the cases when the spectra lie in domains $\tilde{D}_1$ and $\tilde{D}_2$. Obviously, the iteration parameters $\mu_m$ in these cases will be very close to those when the spectra are on diameters $D_1$ and $D_2$. The iteration parameters for $D_2$ will be close to $D_2$ symmetrically with respect to the origin of $S_0$. Without loss of generality, one may assume that segment $D_1$ is on the real axis. Then the equations for ellipses $\tilde{D}_1$ and $\tilde{D}_2$ will have the form

$$\tilde{D}_1(x,y) = \{x,y \mid (x-x_0)^2 + \frac{y^2}{b^2} \leq R^2\},$$

$$(3.16) \qquad \tilde{D}_2(x,y) = \{x,y \mid \frac{(x-x_0)^2}{a^2} + y^2 \leq R^2\},$$

where $x_0 > R$, and $a = a_0, b = b_0, a_0, b_0 \gg 1$.

Consider first the case $n = 2$. Perform a continuous transformation of ellipse $\tilde{D}_1$ first to circle $S_0$ and then to ellipse $\tilde{D}_2$ symmetrically with respect to the origin of $S_0$. Iteration parameters $\mu_1, \mu_2$ will also undergo continuous transformation. When ellipse $\tilde{D}_1$ is transformed to circle $S_0$ the iteration parameters will be a function of $b$, $\mu_1 = \mu_1(b)$, $\mu_2 = \mu_2(b)$, where $b$ varies from $b_0$ to 1. When circle $S_0$ is transformed to ellipse $\tilde{D}_2$ the iteration parameters will be a function of $a$, $\mu_1 = \mu_1(a)$, $\mu_2 = \mu_2(a)$, where $a$ varies from 1 to $a_0$. Taking into account the symmetry, we see that in the process of transformation of a domain, the iteration parameters first merge to a point and then diverge; they merge when the domain is transformed to circle $S_0$ because the direction of the domain transformation changes to a perpendicular one. Then from (3.1), taking into account that the iteration parameters are equal, it follows that the Chebyshev iteration procedure for a circle is identical to simple iterations. Thus we obtain (3.15) for $n = 2$.

Now prove (3.15) for $n = 2^k$, where $k$ is an arbitrary whole number, using the induction. For $k = 1$ the statement was proved. Assume that the statement holds for $k = l$ and prove that it is valid also for $k = l + 1$. Taking into account the symmetry, we see that when ellipse $\tilde{D}_1$ is transformed to circle $S_0$ the iteration parameters $\mu_m, m = 1, ..., 2^{l+1}$ in (2.14) must merge pairwise at $b = 1$. Therefore, from (3.13) we have the following minimax problem for circle $S_0$

$$(3.17) \qquad \rho_{opt} = \min_{\{\mu_m\}} \frac{1}{\left|\prod_{m=1}^{2^l} \mu_m^2\right|} \max_{z \in S_0} \left|\prod_{m=1}^{2^l} (\mu_m - z)^2\right|$$

From (3.17) it follows that this problem is equivalent to the case $n = 2^l$. Thus for all $n$ that are the powers of number 2, relation (3.15) holds. Thus we have

$$(3.18) \qquad \min_{\{\tau_m\}} \max_{z \in S_0} \left|\prod_{m=1}^{n} (1 - \tau_m z)\right| = \prod_{m=1}^{n} \max_{z \in S_0} \left|1 - \frac{1}{\mu_0} z\right| = \left(\frac{R}{x_0}\right)^n$$

The right-hand side of (3.18) specifies, according to Theorem 2.3, the convergence of $n$ successive iterations in GSI.

Now prove *ad absurdum* that (3.15) holds for all $n$. Assume that there exists a number $n_0$, for which (3.15) is not valid. Then there exists iteration parameters $\tau_m^{(n_0)}, m = 1, ..., n_0$, for which the inequality

$$(3.19) \qquad \max_{z \in S_0} \left|\prod_{m=1}^{n_0} (1 - \tau_m^{(n_0)} z)\right| < \left(\frac{R}{x_0}\right)^{n_0}$$

holds. Take a $k_0$ such that $2^{k_0} = n_0 + n_1$, where $n_1 > 0$. Consider the set of iteration parameters

$$\{ \tau_m = \tau_m^{(n_0)}, m = 1, ..., n_0; \ \tau_m = 1/\mu_0, m = n_0 + 1, ..., 2^{k_0} \}.$$

Using (3.19), we have the chain of inequalities

$$\max_{z \in S_0} \left|\prod_{m=1}^{2^{k_0}} (1 - \tau_m z)\right| \leq \max_{z \in S_0} \left|\prod_{m=1}^{n_0} (1 - \tau_m z)\right| \prod_{m=n_0+1}^{2^{k_0}} \max_{z \in S_0} \left|(1 - \frac{1}{\mu_0} z)\right| <$$

$$(3.20) \qquad \left(\frac{R}{x_0}\right)^{n_0} \prod_{m=n_0+1}^{2^{k_0}} \max_{z \in S_0} \left|(1 - \frac{1}{\mu_0} z)\right| = \left(\frac{R}{x_0}\right)^{2^{k_0}}$$

Comparing (3.18) and (3.20), we arrive at contradiction. Thus, if the spectrum of the operator is a circle $S_0$, equality (3.15) holds for all $n$.

In the general case, when the spectrum of the operator lies in the complex plane, an algorithm of determination of iteration parameters is not known and it is necessary to apply numerical methods.

**4. Statement of electromagnetic scattering problems.** Consider the following class of the electromagnetic scattering problems. Assume that in a bounded three-dimensional domain $Q$ the medium is characterized by a permittivity tensor $\hat{\varepsilon}$ (a 3x3 matrix) and its components are functions of the coordinates. Outside domain $Q$ the medium is isotropic and its parameters are constant, i.e. everywhere $\varepsilon = \varepsilon_0 = const$ and $\mu = \mu_0 = const$. It is necessary to determine the electromagnetic field excited by an external field with the time dependence $\exp(-i\omega t)$; both an incident plane wave and current $\mathbf{J}^0$ may be the source of external field. The statement of the corresponding mathematical problem is as follows: find vector-functions $\mathbf{E}$ and $\mathbf{H}$ that satisfy the Maxwell equations

$$(4.1) \qquad \text{curl } \mathbf{H} = -i\omega\hat{\varepsilon}\mathbf{E} + \mathbf{J}^0, \qquad \text{curl } \mathbf{E} = i\omega\mu_0 \mathbf{H}$$

and the radiation conditions at infinity

$$\lim_{r \to \infty}\left[ r\left(\frac{\partial u}{\partial r} - ik_0 u\right)\right] = 0, \quad r = |x| = \sqrt{(x_1^2 + x_2^2 + x_3^2)},$$

where $k_0 = \omega\sqrt{\varepsilon_0\mu_0}$. In (4.1) $\mathbf{J}^0$ is the given current generating external field $\mathbf{E}^0$, $\mathbf{H}^0$ and, in line with the physical essence of the problem, $\text{Im }\varepsilon_0 \geq 0$, $\text{Im }\mu_0 \geq 0$, and $\text{Im }k_0 \geq 0$.

Rewrite equations (4.1) in the equivalent form

$$(4.2) \qquad \text{curl } \mathbf{H} = -i\omega\varepsilon_0 \mathbf{E} + \mathbf{J}, \quad \text{curl } \mathbf{E} = i\omega\mu_0 \mathbf{H}.$$
$$\mathbf{J} = \mathbf{J}^0 + \mathbf{J}^p, \quad \mathbf{J}^p = -i\omega(\hat{\varepsilon} - \varepsilon_0 \hat{I})\mathbf{E}.$$

In (4.2), $\mathbf{J}^p$ is the electric polarization current which does not vanish only in domain $Q$.

We may formally consider (4.2) as the Maxwell equations in a homogeneous medium, assuming that the electromagnetic field is generated by current $\mathbf{J}$. Then a solution to (4.2) satisfying the radiation conditions at infinity can be expressed in terms of vector potential $\mathbf{A}$ using the known formulas [10]

$$\mathbf{A}(x) = \int \mathbf{J}(y)G(R)dy,$$
$$(4.3) \qquad \mathbf{E} = i\omega\mu_0 \mathbf{A} - \frac{1}{i\omega\varepsilon_0}\text{grad div }\mathbf{A}, \quad \mathbf{H} = \text{curl }\mathbf{A}.$$

In (4.3)

$$(4.4) \qquad G(R) = \frac{\exp(ik_0 R)}{4\pi R},$$

where $R = |x - y|$, is the Green's function of the Helmholtz equation. From (4.2)—(4.4) we obtain that the unknown electromagnetic field for the problem under study can be represented in the form

$$\mathbf{E}(x) = \mathbf{E}^0(x) + k_0^2 \int_Q (\hat{\varepsilon}_r(y) - \hat{I})\mathbf{E}(y)G(R)\,dy + \text{grad div}\int_Q (\hat{\varepsilon}_r(y) - \hat{I})\mathbf{E}(y)G(R)\,dy, \; x \in E_3$$
$$(4.5) \qquad \mathbf{H}(x) = \mathbf{H}^0(x) - i\omega\varepsilon_0 \text{ curl}\int_Q (\hat{\varepsilon}_r(y) - \hat{I})\mathbf{E}(y)G(R)\,dy, \; x \in E_3,$$

where $x = (x_1, x_2, x_3)$, $y = (y_1, y_2, y_3)$, and $\hat{\varepsilon}_r = \hat{\varepsilon}/\varepsilon_0$.

Below we will denote the permittivity tensor $\hat{\varepsilon}_r$ by $\hat{\varepsilon}$. Since $\hat{\varepsilon} = \hat{I}$ ($\hat{I}$ is the identity tensor) outside domain $Q$, we can reduce the problem to a volume integrodifferential equation with respect to electric field $\vec{E}$ in domain $Q$

(4.6) $$\mathbf{E}(x) - k_0^2 \int_Q (\hat{\varepsilon} - \hat{I})\mathbf{E}(y)G(R)\,dy - \text{grad div} \int_Q (\hat{\varepsilon} - \hat{I})\mathbf{E}(y)G(R)\,dy = \mathbf{E}^0(x), \ x \in Q.$$

If we find a solution to (4.6) in domain $Q$, we can determine the electromagnetic field outside the domain using integral representations (4.5).

Note that operation *grad div* cannot be applied under the integral sign in (4.6) because double differentiation of $G$ with respect to the coordinates yields the kernel singularity of the order $\sim 1/R^3$ and the resulting integrals diverge. However, outside domain $Q$ operation *grad div* in (4.5) can be applied under the integral sign.

Using the theorem concerning differential properties of weakly singular integrals [7] we can reduce equation (4.6) to a VSIE [9]

(4.7) $$\mathbf{E}(x) + \frac{1}{3}(\hat{\varepsilon}(x) - \hat{I})\mathbf{E}(x) - p.v.\int_Q ((\hat{\varepsilon}(y) - \hat{I})\mathbf{E}(y), \text{grad})\,\text{grad}\,G(R)\,dy -$$
$$k_0^2 \int_Q (\hat{\varepsilon}(y) - \hat{I})\mathbf{E}(y)G(R)\,dy = \mathbf{E}^0(x), \ x \in Q.$$

Here symbol $(*,*)$ denotes the inner product of vectors and $p.v.\int$ a singular integral, i.e. an integral over a domain where an infinitesimal ball centered in a vicinity of the point $x = y$ is extracted.

Note that equation (4.7) describes the scattering problems with minimal restrictions imposed on parameters of the medium; namely, it is assumed that the permittivity is a bounded function of the coordinates. However, when the spectrum of the operator is studied one has to impose certain restrictions on the permittivity tensor required by the fact that we apply in our proofs the results of the theory of singular equations. Taking into account these restrictions, we will assume in what follows that the tensor components are Hölder-continuous functions of the coordinates. We use in this connection the following

DEFINITION. *$u$ is a Hölder-continuous function in domain $D$ if the inequality*

$$|u(y) - u(x)| \leq C|x - y|^\delta, \quad C = \text{const}, 1 \geq \delta > 0,$$

*holds for any $x, y \in D$.*

**5. Spectrum of the integral operator.** The spectrum of operator $\hat{A}$ is a set $\lambda$ of the points on the complex plane $Z$ such that the operator $(\hat{A} - \lambda \hat{I})$ is not invertible everywhere in the Hilbert space $H$. The points of $\lambda$ at which $(\hat{A} - \lambda \hat{I})$ is not a Noether operator belongs to the continuous part of the spectrum (essential spectrum) of $\hat{A}$. The points of $\lambda$ at which there exists a nontrivial solution $u$ of the homogeneous equation $\hat{A}u - \lambda u = 0$ belongs to the discrete part of the spectrum of $\hat{A}$ [3].

First we choose an appropriate Hilbert space. The integrals of the squared moduli of the characteristics of the electromagnetic field enter the energy conservation law (the Poynting theorem). Therefore, the space of square-integrable functions is the most 'physical' space as far as the analysis of integral equations associated with electromagnetic scattering problems are concerned. Below, we will use Hilbert space $\mathbf{L}_2$ with the inner product

$$(\mathbf{U}, \mathbf{V}) = \int \mathbf{U}(x) \mathbf{V}^*(x)\,dx.$$

The following statement holds [9]

THEOREM 5.1. *Assume that the Cartesian components of tensor $\hat{\varepsilon}(x)$ are Hölder-continuous functions everywhere in Euclidean space $E_3$. Then the VSIE operator in (4.7) is a Noether operator in $\mathbf{L}_2(Q)$ if and only if the following condition holds*

(5.1) $$\sum_{n,m=1}^{3} \varepsilon_{nm}(x) \beta_n \beta_m \neq 0$$

for all $x \in Q$ and real numbers $\beta_1, \beta_2, \beta_3$ such that $\beta_1^2 + \beta_2^2 + \beta_3^2 = 1$.

For an isotropic medium, equation (5.1) takes the form

$$\varepsilon(x) \neq 0, \ x \in Q.$$

Write integral equation (4.7) in a symbolic form

(5.2) $$\hat{A}u \equiv u - \hat{S}((\hat{\varepsilon} - I)u) = f,$$

where operator $\hat{S}$ is given in (4.7). Obviously,

(5.3) $$\hat{A} - \lambda \hat{I} = (1 - \lambda)\left[\hat{I} - \hat{S}\left(\frac{\hat{\varepsilon} - \lambda I}{1 - \lambda} - I\right)\right].$$

Introduce a tensor-function

(5.4) $$\hat{\varepsilon}^+(\lambda, x) = \frac{\hat{\varepsilon}(x) - \lambda \hat{I}}{1 - \lambda}.$$

Substituting (5.4) into (5.3), we obtain

(5.5) $$\hat{A} - \lambda \hat{I} = (1 - \lambda)(\hat{I} - \hat{S}(\hat{\varepsilon}^+ - I)).$$

Comparing (5.5) and (5.2), taking into account (5.4) and the relation $\beta_1^2 + \beta_2^2 + \beta_3^2 = 1$, and applying Theorem 5.1 we see that the set $\sigma_1$ on complex plane $Z$ defined from the conditions

(5.6) $$\lambda = \sum_{n,m=1}^{3} \varepsilon_{nm}(x) \beta_n \beta_m, \ x \in Q, \ \beta_1^2 + \beta_2^2 + \beta_3^2 = 1$$

belongs to the continuous part of the spectrum of the operator of equation (4.7). The latter condition implies that the point $\lambda = 1$ belongs to $\sigma_1$; in fact, by virtue of the Hölder-continuity, the permittivity tensor becomes a scalar quantity on the boundary of domain $Q$, i.e. $\varepsilon_{nm} = \delta_{nm}$.

For an isotropic medium, we have the following formula for the points of the continuous spectrum

(5.7) $$\lambda = \varepsilon(x), \ x \in Q.$$

Denote by $\sigma$ the minimal simply connected set on complex plane $Z$ containing $\sigma_1$. Consider a simply connected set $\sigma^+ = Z \setminus \sigma$. Then taking into consideration the foregoing analysis, we can state that $(\hat{A} - \lambda \hat{I})$ is a Noether operator if $\lambda \in \sigma^+$. It is known [5] that in every connected component of the domain where an operator is of Noether type, the operators $(\hat{A} - \lambda \hat{I})$ have the same index. If $|\lambda| > \|\hat{A}\|$, then operator $(\hat{A} - \lambda \hat{I})$ has a bounded inverse defined everywhere; therefore the index of the operator is zero. Since domain $\sigma$ is bounded there exist points $\lambda_0 \in \sigma^+$ such that $|\lambda_0| > \|\hat{A}\|$; consequently, the index of (Noether) operators $(\hat{A} - \lambda \hat{I})$ equals zero at $\lambda \in \sigma^+$ and they are Fredholm operators. We have proved the following

THEOREM 5.2. *The set $\sigma_1$ on the complex plane defined by (5.6) belongs to the continuous spectrum of the operator of integral equation (4.7). In addition, $(\hat{A} - \lambda \hat{I})$ is a Fredholm operator if $\lambda \in \sigma^+ = Z \setminus \sigma$, where $\sigma$ is the minima simply connected set containing $\sigma_1$.*

Note that domain $\sigma_1$ of the continuous spectrum is governed solely by the permittivity values and is independent of geometrical properties of domain $Q$.

From Theorem 5.2 it follows that points $\lambda \in \sigma^+$ belong either to the resolvent set of operator $\hat{A}$ or to its discrete spectrum. Generally, it is not possible to describe sufficiently accurate the domain of localization of the discrete spectrum of operator (4.7). However, we can do it in one very important particular case.

Below we consider the problem of the low-frequency scattering of electromagnetic waves when diameter $D$ of domain $Q$ is much less than the wavelength, i.e. $D << \dfrac{2\pi}{k_0}$.

Equation (4.7) may be also considered in the static case when the wavenumber $k_0 = 0$. This circumstance shows an essential difference between three-dimensional and two-dimensional problems, because in the latter case a transformation to the static case of stationary integral equations cannot be performed. From (4.7) we obtain

(5.8) $\quad (\hat{A}(k_0) - \hat{A}(0))\mathbf{E} = -k_0^2 \int_Q (\hat{\varepsilon} - \hat{I})\mathbf{E}(y) G(R) dy - \int_Q ((\hat{\varepsilon} - \hat{I})\mathbf{E}(y), \text{grad}) \, \text{grad} \, G_0(R) dy$,

where $\hat{A}(k_0)$ and $\hat{A}(0)$ are the operators of integral equations for the stationary and static cases and

(5.9) $\quad\quad\quad\quad\quad\quad G_0(R) = \dfrac{\exp(ik_0 R) - 1}{4\pi R}$.

The second integral operator in (5.8) is not a singular one because due to (5.9) its kernel has no singularity at $x = y$ and is a smooth function of coordinates. Therefore, from (5.8) it follows that

$$\lim_{k_0 \to 0} \|\hat{A}(k_0) - \hat{A}(0)\| = 0.$$

Then, using the known result [5] of functional analysis concerning convergence of the spectra of operators we have the following

LEMMA 5.1. *The spectrum of low-frequency integral operator $\hat{A}(k_0)$ converges to the spectrum of static integral operator $\hat{A}(0)$ when $k_0 \to 0$.*

In the static case, integrodifferential equation (4.6) which is equivalent to singular integral equation (4.7) can be written in the form

(5.10) $\quad\quad\quad\quad \mathbf{E}(x) - \text{grad div} \int_Q (\hat{\varepsilon}(y) - I)\mathbf{E}(y)(1/4\pi R)) \, dy = \mathbf{E}^0(x)$.

The solution of homogeneous equation (5.10) at $\mathbf{E}^0 = 0$ satisfies the differential equations
(5.11) $\quad\quad\quad\quad\quad\quad \text{curl } \mathbf{E} = 0, \quad \text{div}(\hat{\varepsilon} \mathbf{E}) = 0$.

The first equation in (4.11) follows from the identity *rot grad = 0*. The second equation follows from the identities *grad div = curl curl + $\Delta$* and *div curl = 0* and the differential equation $\Delta \mathbf{A} = -\mathbf{J}$ which is satisfied by the volume potential $\mathbf{A}(x) = \int \mathbf{J}(y)(1/4\pi R) dy$ [10].

From the first equation in (5.11) we have $\mathbf{E} = \text{grad } \varphi$. Consequently, equations (5.11) are reduced to a second-order differential equation with respect to a scalar function $\varphi$

(5.12) $\quad\quad\quad\quad\quad\quad \text{div}(\hat{\varepsilon} \, \text{grad } \varphi) = 0$.

Let $\psi$ be a differentiable function defined everywhere. Then the following obvious identity holds

(5.13) $$\operatorname{div}(\psi\hat{\varepsilon}\operatorname{grad}\varphi) = \psi\operatorname{div}(\hat{\varepsilon}\operatorname{grad}\varphi) + \operatorname{grad}\psi \bullet \hat{\varepsilon}\operatorname{grad}\varphi.$$

Set $\psi = \bar{\varphi}$. Then integrating (5.13) over the whole space and taking into account (5.12) and the divergence theorem, we obtain an integral relationship

(5.14) $$\int (\hat{\varepsilon}\operatorname{grad}\varphi, \operatorname{grad}\varphi)\, dv = \lim_{R\to\infty}\int_{S_R} \bar{\varphi}\frac{\partial\varphi}{\partial n}\, dS,$$

where $S_R$ is a sphere of radius $R$ centered at the origin and $n$ is the normal vector to the sphere. Since $\varphi$ is a harmonic function outside domain $Q$, the integrand $\bar{\varphi}\,\partial\varphi/\partial n$ decays at infinity not slower than $R^{-3}$. Thus, the limit on the right-hand side of (5.14) equals zero and every solution of homogeneous equation (5.12) and consequently of homogeneous integral equation (5.10) with $\mathbf{E} = \operatorname{grad}\varphi$ satisfies the integral relationship

(5.15) $$\int (\hat{\varepsilon}\operatorname{grad}\varphi, \operatorname{grad}\varphi)\, dv = 0.$$

From (5.5), it is clear that $\lambda$ will be a point of the discrete spectrum of the operator of equation (5.10) if there exists a nonzero solution $\varphi(\lambda, x)$ of equation (5.12) with the permittivity function (5.4). Then, from (5.15) and (5.4) it follows that the corresponding $\lambda$ can be determined from the formula

(5.16) $$\lambda = \frac{\int (\hat{\varepsilon}(x)\operatorname{grad}\varphi(\lambda, x), \operatorname{grad}\varphi(\lambda, x))\, dx}{\int |\operatorname{grad}\varphi(\lambda, x)|^2\, dx}.$$

Obviously it is not possible to determine functions $\varphi(\lambda, x)$. However one can use (5.16) to find a domain of localization of the points of the discrete spectrum on complex plane Z.

Consider first the isotropic case. From (5.16) we obtain that of the points of the discrete spectrum are given by

(5.17) $$\lambda = \frac{\int \varepsilon |\operatorname{grad}\varphi|^2\, dx}{\int |\operatorname{grad}\varphi|^2\, dx}.$$

Expression (5.17) is closely related to the formula specifying the center of mass of a plane figure which can be located only inside a convex envelope of the figure. One can show in a similar manner that the values of $\lambda$ is situated only inside a convex envelope of the domain of permittivity $\varepsilon(x)$, $x \in Q$, i.e. inside a convex envelope of $\sigma_1$. Denote this set by $\sigma_0$. Obviously, the inclusion $\sigma_0 \supseteq \sigma$ is valid, where $\sigma$ is determined according to Theorem 5.2. Therefore, domain $Z \setminus \sigma_0$ may contain only the points of the discrete spectruma which are not situated there according to the analysis performed above. We have proved the following

THEOREM 5.3. *In the case of isotropic medium the whole spectrum of the operator of equation (5.10) and therefore of singular operator (4.7) at $k_0 = 0$ is situated inside a convex envelope of the set given by formula (5.7).*

For the anisotropic case the situation is more complicated because expression (5.16) cannot be treated in the same manner as in the isotropic case. Namely, write tensor-function $\hat{\varepsilon}(x)$ as

(5.18) $$\hat{\varepsilon}(x) = \hat{\delta}_1(x) + i\,\hat{\delta}_2(x), \quad \hat{\delta}_1(x) = \frac{\hat{\varepsilon}(x) + \hat{\varepsilon}^*(x)}{2}, \quad \hat{\delta}_2(x) = \frac{\hat{\varepsilon}(x) - \hat{\varepsilon}^*(x)}{2i}.$$

$\hat{\delta}_1$ and $\hat{\delta}_2$ are Hermitian tensors for all $x \in Q$ so that their eigenvalues are real numbers. Denote by $a_{\min}^{(1)}(x), a_{\max}^{(1)}(x)$ and $a_{\min}^{(2)}(x), a_{\max}^{(2)}(x)$ the minimal and maximal eigenvalues of Hermitian tensors $\hat{\delta}_1(x)$ and $\hat{\delta}_2(x)$.

From (5.16) and (5.18) it follows

(5.19) $$\lambda = \lambda_1 + i\lambda_2, \quad \lambda_k = \frac{\int (\hat{\delta}_k(x)\,\text{grad}\,\varphi(\lambda,x), \text{grad}\,\varphi(\lambda,x))\,dx}{\int |\text{grad}\,\varphi(\lambda,x)|^2\,dx}, \quad k=1,2.$$

Since $(\hat{\delta}_k \vec{V}, \vec{V})$, $k=1,2$ are real quantities for any three-dimensional complex vectors $\vec{V}$, then $\lambda_1$ and $\lambda_2$ in (5.19) are also real quantities. Next, for all $x$ we have

$$a_{\min}^{(k)}(x)|\text{grad}\,\varphi(\lambda,x)|^2 \leq (\hat{\delta}_k(x)\,\text{grad}\,\varphi(\lambda,x), \text{grad}\,\varphi(\lambda,x)) \leq a_{\max}^{(k)}(x)|\text{grad}\,\varphi(\lambda,x)|^2, \quad k=1,2.$$
(5.20)

Now, substituting (5.20) into (5.19), we obtain the following

THEOREM 5.4. *The whole spectrum of the operator of equation (5.10) and therefore of singular operator (5.7) at $k_0 = 0$ is situated inside a rectangle whose sides are parallel to the axes of the complex plane and the left lower and right upper vertices have the coordinates* $(A_{\min}^{(1)}, A_{\min}^{(2)})$, $(A_{\max}^{(1)}, A_{\max}^{(2)})$, *where*

$$A_{\min}^{(n)} = \min a_{\min}^{(n)}(x), \quad A_{\max}^{(n)} = \max a_{\max}^{(n)}(x), \, x \in Q, n = 1,2.$$

Note that Theorem 5.4 gives less exact information about the domain of localization of the spectrum as compared with the isotropic case. From Theorems 5.3 and 5.4 it follows that if the medium is lossless, then the whole spectrum of the operator lies close to a segment of the real axis of the complex plane.

Lemma 5.1 establishes a possibility of using the obtained information about the spectrum for solving the problem of low-frequency electromagnetic scattering with the help of VSIEs and stationary iteration methods set forth in Sections 2 and 3.

All the results in this Section are proved under the assumption that the parameters of the medium are given as Hölder-continuous functions. However, these restrictions are not essential from the practical viewpoint. In fact, one really important family of the media where the parameters are not described by Hölder-continuous functions is a class of piecewise continuous permittivity functions. A domain occupied by the scatterer with a constant value of the permittivity constitutes the simplest example. Note in this respect that in the vicinity of the interface the permittivity can be represented as a Hölder-continuous function which undergoes a finite break along an infinitesimal segment. Then, all the results can be applied in the case of piecewise continuous permittivity functions (piecewise smooth media).

Explain the above conclusion by using the following example. Let domain $Q$ in equation (4.7) be a ball and the permittivity function (in the isotropic case) have the following form in spherical coordinates

(5.21) $$\varepsilon(r) = \begin{cases} \varepsilon_2, & d_2 \geq r \geq 0 \\ \varepsilon_2 + (\varepsilon_1 - \varepsilon_2)\dfrac{r - d_2}{d_1 - d_2}, & d_1 \geq r \geq d_2 \\ \varepsilon_1 + (1 - \varepsilon_1)\dfrac{r - d_1}{R - d_1}, & R \geq r \geq d_1 \end{cases}$$

In (5.21) $R$ is the radius of the ball and $R > d_1 > d_2 > 0$.

The solid line in Fig. 5.1 shows the continuous part of the spectrum while the whole spectrum of the operator associated with a low-frequency scattering problem lies inside a triangle.

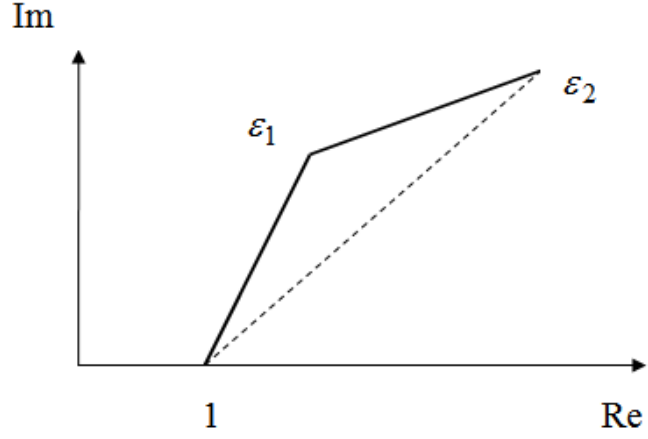

**Figure 5.1:** Spectrum of the integral operator.

When $d_1$ tends to $R$ the continuous part of the spectrum remains unchanged while its discrete part (in the low-frequency case) undergoes variation but remains inside the triangle.

**6. Numerical experiments.** Let us illustrate theoretical findings of this study by the results of computations obtained for a representative problem: low-frequency scattering of a plane electromagnetic wave from an inhomogeneous isotropic dielectric ball; the parameters of the ball taken according to (5.21) are $\varepsilon_2 = 3+i$, $\varepsilon_1 = 2+2i$, $d_2 = R/2$, $d_1 = 2R/3$. Consider the spectra of the VSIE operators. There is a correspondence between the spectrum of the integral operator and its discrete (matrix) analogue [5]. Of course, a matrix has only eigenvalues and has no continuous spectrum; therefore we cannot see continuous (bold) lines as in Fig. 5.1. However the eigenvalues will thicken along the lines of continuous spectrum.

Dicretizing VSIE (4.7), we obtain a SLAE of dimension $N$, where $N$ is determined from the conditions of approximation of the integral operator. Figure 6.1 shows the spectra of the SLAE matrices calculated for three values of the radius of the ball (5.21)

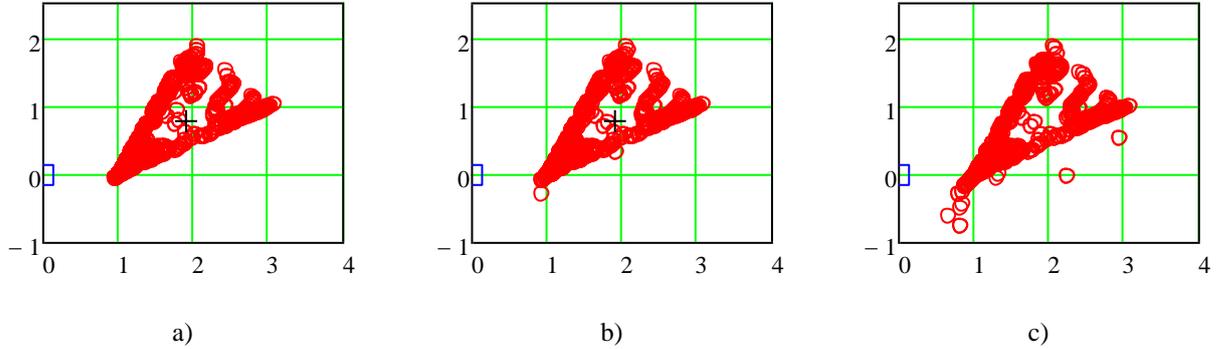

Figure 6.1. Spectra for the ball radii $R = \lambda/20$ (a), $R = \lambda/6$ (b), and $R = \lambda/2$ (c)

One can see from Fig. 6.1 that for all values of the ball dimensions the spectrum of the operator contains the lines corresponding to the VSIE continuous spectrum. The size of the ball in Fig. 6.1 (a) is within a low-frequency range and the whole spectrum of the operator is contained in the triangle with the vertexes (1,0), (3,1), (2,2). In Fig. 6.1 (b) one can see the points of the spectrum of the operator appearing outside the triangle. Fig. 6.1 (c) shows a resonance range and the convex envelope of the spectrum differs considerably from the initial triangle.

The calculated spectra of the SLAE matrices for a homogeneous isotropic dielectric ball of the radius $R = \lambda/30$ are presented in Figs. 6.2 and 6.3.

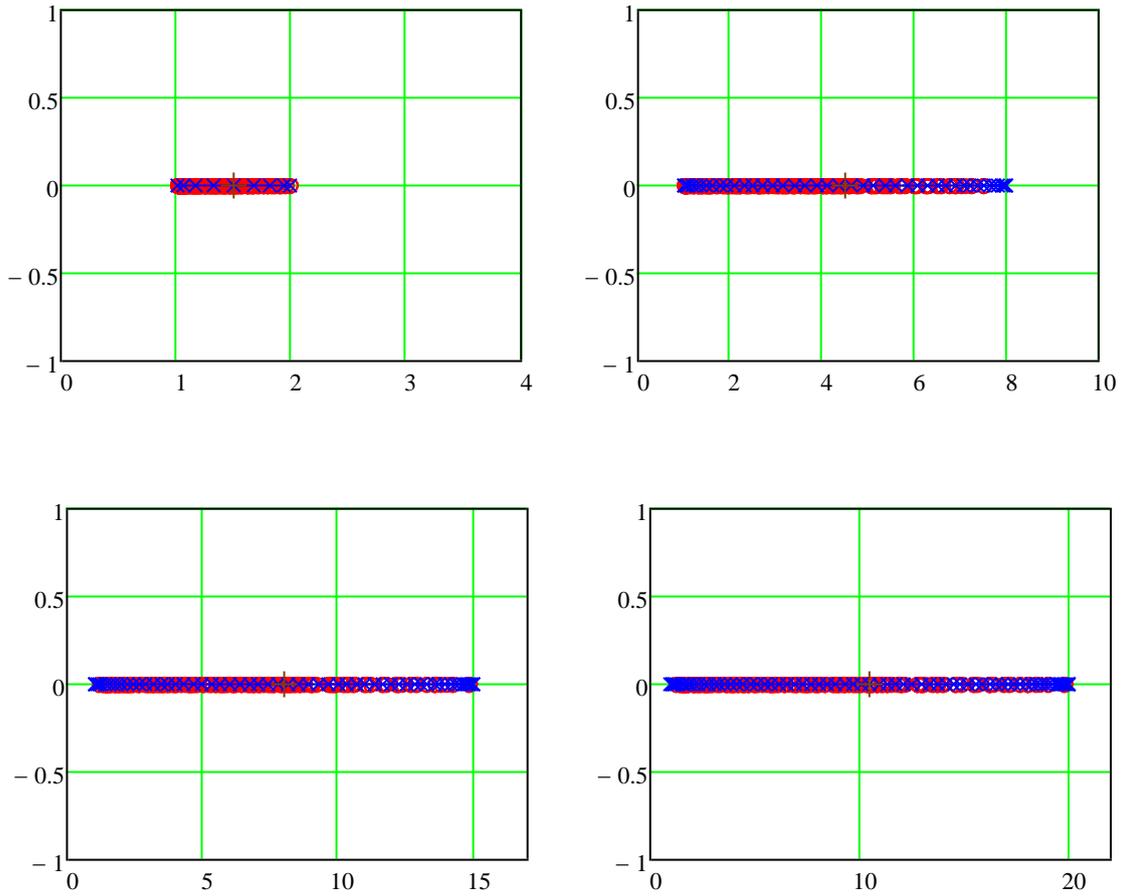

Figure 6.2. Spectra of the SLAE matrices for a homogeneous dielectric ball of the radius $R = \lambda/30$ calculated for the permittivity values $\varepsilon = 2, 8, 15, 20$.

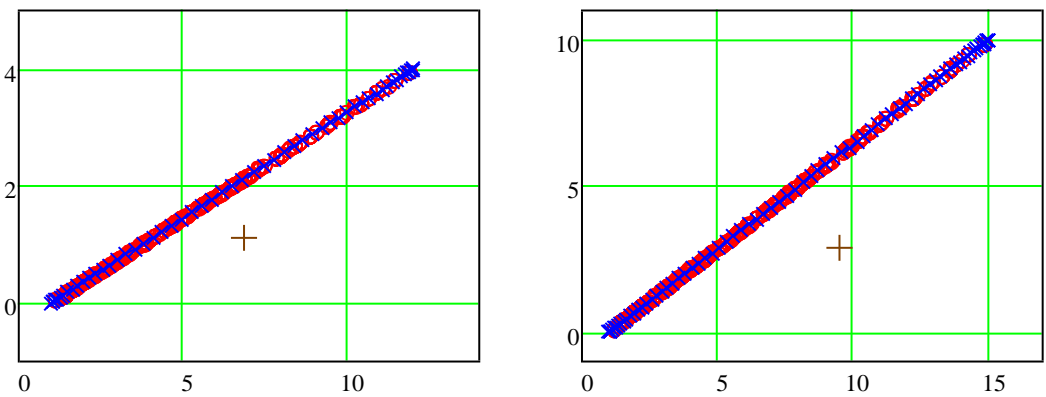

Figure 6.3. Spectra of the SLAE matrices for a homogeneous dielectric ball of the radius $R = \lambda/30$ calculated for the permittivity values $\varepsilon = 12 + 4i, \; 15 + 10i$.

It can be seen that the spectra of the operators lie in the complex plane inside the segment $[1, \varepsilon]$ which agrees with the theory developed in this study. Note that for lossless bodies the spectrum is situated on the real axis.

We also calculated the spectra of the SLAE matrix in the quasi-static range for a homogeneous isotropic dielectric cube having the same permittivities. The obtained graphs coincide with those shown in Figs. 6.2 and 6.3, which confirms, in line with the results of this work, that the spectrum is independent of the shape of the scatterer.

The spectra of the SLAE matrices calculated for an inhomogeneous dielectric ball of the radius $R = \lambda/30$, $d_2 = R/2$, $d_1 = 2R/3$, and permittivity given by (5.21) are shown in. Fig. 6.4.

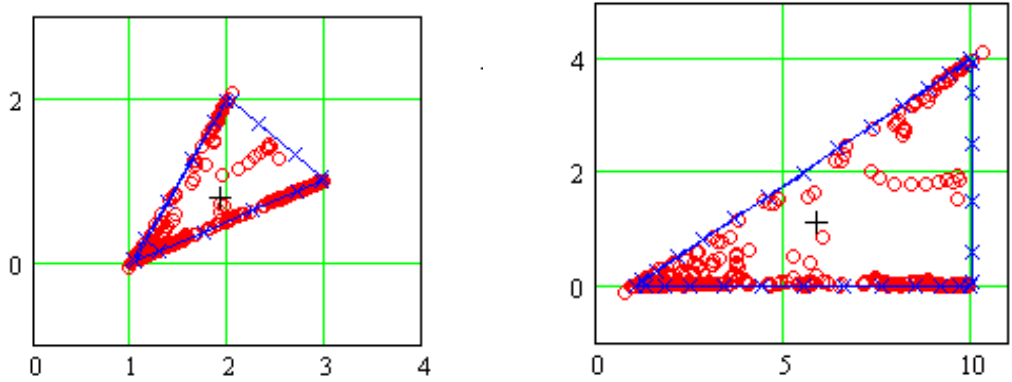

Figure 6.4. Spectra of the SLAE matrices for an inhomogeneous dielectric ball of the radius $R = \lambda/30$ and permittivity given by (5.21) with $\varepsilon_2, \varepsilon_1 =$(2+2i, 3+i), (10+4i, 10).

For other $d_2, d_1$ the complex spectra of the operators determined for the same values of $\varepsilon_2, \varepsilon_1$ are always contained inside the triangle with the vertexes (1,0), $\varepsilon_2, \varepsilon_1$.

Mutual location of the spectra of the VSIE operators and SLAE matrix eigenvalues calculated in [2, 8] also demonstrate complete agreement with the theoretical results presented in this work.

Compare now the convergence of GSI, GCI and GMRES. As a representative problem we choose the scattering of a plane electromagnetic wave from a dielectric ball with the radius $R = \lambda/30$. Table 6.1 gives the number of iterations $L$ required for attaining the relative computational error

$$\delta = \frac{\|A\vec{E} - \vec{E}^0\|}{\|\vec{E}^0\|},$$

where $\hat{A}$ is a discrete analogue of the integral operator and $\vec{E}$ is the calculated approximate solution. Computations were performed until $\delta \leq 10^{-5}$.

Table 6.1.

|  | GSI | GCI, n = 5 | GCI, n = 10 | GMRES, n = 2 | GMRES, n = 5 | GMRES, n = 10 |
|---|---|---|---|---|---|---|
| $\varepsilon = 2$ | L=10 | L=10 | L=10 | L=8 | L=10 | L=10 |
| $\varepsilon = 8$ | L=55 | L=25 | L=30 | L=32 | L=25 | L=30 |
| $\varepsilon = 15$ | L=155 | L=80 | L=80 | L=3044 | L=1175 | L=680 |
| $\varepsilon = 20$ | L=723 | L=240 | L=210 | L=4002 | L=1250 | L=830 |
| $\varepsilon = 12+4i$ | L=55 | L=30 | L=30 | L=84 | L=50 | L=40 |
| $\varepsilon = 15+10i$ | L=76 | L=40 | L=30 | L=144 | L=75 | L=60 |
| $\varepsilon_2 = 2+2i$, $\varepsilon_1 = 3+i$ | L=16 | L=15 | L=10 | L=26 | L=25 | L=20 |

For GCI we present the number of iterations $L$ for two values of parameter $n$ entering iteration procedure (3.1), (3.12). For GMRES we consider three values of the dimension $n$ of the Krylov subspace. The first six rows of Table 6.1 contain the data for the scattering from a lossless and then lossy homogeneous dielectric ball and the last two rows for an inhomogeneous dielectric ball with the permittivity given by (5.21) and $d_2 = R/2$, $d_1 = 2R/3$. For GCI the optimal iteration parameter was determined in all cases according to the algorithm of Section 2. For the lossless homogeneous dielectric ball the GCI iteration parameters were calculated by formula (3.14) because in this case the spectrum of the operator is situated close to the real axis. For the lossy homogeneous dielectric ball the spectrum of the operator is situated on the complex plane close to a segment with the endpoint (1,0); the exact values of the GCI iteration parameters (3.12) are not known. However, their good approximations can be determined by the following algorithm: formula (3.14) is used for calculating the iteration parameters for a segment $S$ of the real axis equally long with complex segment $(1, \varepsilon)$, then segment $S$ is turned around point (1,0) until it merges with complex segment $(1, \varepsilon)$ so that finally the sought-for iteration parameters are located on the complex segment. For an inhomogeneous dielectric ball the spectrum is located inside the triangle for which the exact values of the GCI iteration parameters are not known either; in this case we determine approximate values of the optimal iteration parameters that lie on the sides of the triangle using the algorithm described above for a segment. The data of Table 6.1 clearly shows that the GMRES convergence rate is close to that of GSI and GCI only for moderate values of the permittivity.

We see that as the permittivity increases, the stationary iteration methods with optimal iteration parameters determined according to the algorithms developed in this work turn out to be *much more efficient*, sometimes by an order of magnitude (see e.g. the fourth row of Table 6.1), at least for the low-frequency scattering problems. In addition, the memory volume required for implementation of the stationary algorithms (the dimension of the Krylov subspace) is approximately $n$ times smaller than that for GMRES.

**7. Discussion and conclusions.** The theory of GSI (Section 2) is developed in its complete form. In fact, for an arbitrary domain of the spectrum localization on the complex plane we have constructed a finite algorithm of finding the optimal iteration parameter.

In Section 3 we have described GCI applied for solving linear equations with nonselfadjoint operators. We have obtained sufficient conditions providing the convergence of iterations imposed on the domain of localization of the spectrum on the complex plane. A minimax problem for the determination of optimal complex iteration parameters has been formulated. Note that the memory volume required for the GCI implementation does not depend on the number $n$ of layers in iteration procedure (3.1) and is governed, as well as for simple iterations, by the relationship $M_{ITER} \approx N$, where $N$ is the SLAE dimension.

We have obtained the following important results: if the spectrum of the operator is a circle $S_0$ on the complex plane, the problem degenerates. It means that all iteration parameters merge and become equal to the GSI optimal iteration parameter. Also it is reasonable to take small $n \approx 3-5$ when GMRES is used. Greater dimensions lead mainly to increasing the memory volume rather than to significant improvement of convergence. Our numerical investigations including the results of Section 6 confirmed this conclusion several times.

In Section 4 we set forth VSIEs associated with the problems of electromagnetic wave scattering by dielectric structures. In Section 5 we study the spectrum of integral equations. In the general case and for the resonance wavelength range, we describe explicitly the continuous part of the spectrum of the operator on the complex plane. We would like to strengthen that this part of the spectrum depends solely on the permittivity and is independent of the shape of the scatterer which is confirmed by the numerical results of Session 6. Note that it is not possible to obtain sufficiently exact information on the location of the discrete spectrum in the general case, and in particular, for the scattering problems in the resonance wavelength range. For the low-frequency scattering problems we determine highly accurate the localization of the whole spectrum including its discrete part. Theoretical estimates of the spectrum location on the complex plane agree completely with the calculated eigenvalues of the SLAE matrix resulting from VSIE discretization and presented in Section 6. If the medium is lossless (or losses are negligible) the formulas obtained in this study show that the spectrum of the operator lies close to a segment of the real axis and one can use optimal iteration parameters given by (3.14).

Let us formulate the main conclusion of this study: from the viewpoint of computational resources, the considered stationary iteration methods are the most efficient techniques for the numerical solution of VSIEs describing three-dimensional low-frequency electromagnetic wave scattering by dielectric structures.

**Acknowledgements.** This work was partly supported by the Visby research program of the Swedish Institute.

# REFERENCES

[1] N.V. Budko, A.B. Samokhin, and A.A. Samokhin. *A Generalized overrelaxation method for solving singular volume integral equations in low-frequency scattering problems,* Differential Equations, 41 (2005) 1262-1266.

[2] N.V. Budko, A.B. Samokhin. *Spectrum of the volume integral operator of electromagnetic Scattering,* SIAM J. Sci. Comput., 28 (2006) 682-700.

[3] N. Danford. J. Schwartz. *Linear operators, Spectral operators*, Wiley, 1988.

[4] M.H. Gutknecht. *Lanczos-type solvers for nonsymmetric linear system of equations*. Swiss Center for scientific computing, Zurich, 1997.

[5] L.V. Kantorovich, G.P. Akilov. *Functional analysis*, Pergamon Press, Oxford, 1982.

[6] R.E. Kleinman, P.M. van den Berg. *Iterative methods for radio-wave problems*, The Review of Radio Science, Oxford University Press, (1993), 54-74.

[7]. K. Kobayashi, Y. Shestopalov, Y. Smirnov. *Investigation of electromagnetic diffraction by a dielectric body in a waveguide using the method of volume singular integral equation,* SIAM J. Appl. Maths., 70 (2009) 969-983

[8] G.I. Marchuk. *Methods of numerical mathematics*, Springer, Berlin, 1982.

[9] S.G. Mikhlin. S. Prössdorf. *Singular integral operators,* Springer, Berlin, 1986.

[10] J. Rahola, *On the eigenvalues of the volume integral operator of electromagnetic scattering,* SIAM J. Science Comput., 21 (2000) 1740-1754.

[11] Y. Saad, M.N. Schultz. GMREZ: *A generalized minimal residual algorithm for solving nonsymmetric linear systems*, SIAM J. Sci. Stat. Comput., 7 (1986) 856-869.

[12] A.B.Samokhin. *Integral equations and iteration methods in electromagnetic scattering*, VSP, Utrecht, 2001.

[13] T.K. Sarkar (Ed.) *Application of conjugate gradient method in electromagnetics and signal analysis*, Elsevier, New York, 1991.

[14]. Y. Shestopalov, Y. Smirnov. *Existence and uniqueness of solution to the inverse problem of complex permittivity reconstruction of a dielectric body in a waveguide* Inverse Problems, 26 (2010) 105002

[15]. Y. Shestopalov, Y. Smirnov. *Determination of permittivity of an inhomogeneous dielectric body in a waveguide* Inverse Problem, 27 (2011) 095010

[16] J.A. Stratton. *Electromagnetic theory*, McGraw-Hill, New York, 1941.